\documentclass[reqno]{article}

\usepackage{amsmath}
\usepackage{amssymb}

\begin{document}

\title{ Multiresolution in the Bergman space }
\author{Margit Pap,\\
Marie Curie fellow, NuHAG\\
Faculty of Mathematics,
University of Vienna\\
Alserbachstra\ss e 23,
A-1090 Wien, Austria\\
and\\
University of P\'ecs,\\
 Ifj\'us\'ag \'utja 6,
 7634 P\'ecs, Hungary}
\date{}
\maketitle \centerline{ papm@gamma.ttk.pte.hu, margit.pap@univie.ac.at}
\vskip3mm

 \vskip3mm

% \centerline{This research was supported by
%the Hungarian National Science Foundation}

%\centerline{T047128 OTKA and NKFP-2/020/2004.}

\smallskip

%\keywords
%
%\endkeywords
%\subjclass
%
%\endsubjclass

\abstract
{In this paper we give  a   multiresolution construction in Bergman space. The successful application of rational orthogonal bases needs a priori knowledge of the poles of the transfer function that may cause a drawback of the method.   We give a set of poles  and using them we will generate a multiresolution in $A^2$. We  give sufficient condition for this set to be  sampling sequence for the Bergman space. The construction is an analogy with the discrete affine wavelets, and in fact is the discretization of the continuous voice transform generated by a representation of the Blaschke group over the Bergman space. The constructed discretization scheme gives opportunity of practical realization of hyperbolic wavelet representation of signals belonging to the Bergman space    if we can measure the values of the transfer function on a given set of points inside the unit disc. Convergence properties of the hyperbolic wavelet representation will be studied.}



\vskip4mm \centerline{\bf 1. Introduction}
The plan of this paper is as follows. First we  introduce  a discrete subset  of the Blaschke group.  We will give sufficient conditions for this discrete subset to be  sampling set for the Bergman space. We present some basic results connected to the Bergman space, we give the  definition
of the voice transform  generated by a representation of the Blaschke group on $A^2$.
   Using the discrete subset of the Blaschke group
  we construct a  multiresolution decomposition in $A^2$.  First the different resolution spaces will be defined
   using  nonorthogonal basis which shows the analogy between the discrete hyperbolic wavelets in $A^2$ and
    the discrete affine wavelets in $L^2(\Bbb R)$. Applying  the Gram-Schmidt orthogonalization we consider the
     rational orthogonal basis on the $n$-th multiresolution level $V_n$. This system is the analogue of
     the Mamquist-Takenaka system in the Hardy spaces, possesses similar properties and is connected to the contractive zero divisors of a finite set in Bergman space.  We prove that the projection operator $ P_nf(z)$ on the resolution level $V_n$ 
  is convergent in $A^2$ norm to $f$,  is interpolation operator on the
  set the $\bigcup_{k=0}^n{\cal{ A}}_k$, where ${\cal{ A}}_k$  is defined by (2.7) with minimal norm and $ P_nf(z)\to  f(z)$
  uniformly on every compact subset of the unit disc.

\vskip4mm

\centerline  {\bf 1.2. The Blaschke group}

\vskip4mm

Let us denote by
$$
B_{{\bf a}}(z):=\epsilon\frac{z-b}{1-\bar b z}\ \ (z\in \mathbb C,
{\bf a}=(b,\epsilon)\in \mathbb B:=\mathbb D\times \mathbb T,
\overline{b}z\ne 1) \leqno (1.1)
$$
 the so called  {\it Blaschke functions}, where
$$\mathbb D_+:=\mathbb D:=\{z\in\mathbb C:|z|<1\},\ \ \mathbb T:=\{z\in\mathbb C:|z|=1\}. \leqno ( 1.2)
$$
If ${\bf a} \in \mathbb B$, then $B_{\bf a}$ is an 1-1 map on $\mathbb T
$,and \,\,$\mathbb D $ respectively. The
restrictions of the Blaschke functions on the set $\mathbb D$ or on
$\mathbb T$ with the operation $(B_{{\bf a}_1}\circ
B_{{\bf a}_2})(z):=B_{{\bf a}_1}(B_{{\bf a}_2}(z))$ form a group. In the set of the
parameters $\mathbb B:=\mathbb D\times\mathbb T$ let us define the
operation induced by the function composition in the following way:
$B_{{\bf a}_1}\circ B_{{\bf a}_2}=B_{{\bf a}_1\circ {\bf a}_2}$. The group $(\mathbb B,
\circ)$ will be isomorphic with the group $(\{B_{\bf a}, {\bf a}\in \mathbb B
\}, \circ)$.
 If we use the notations ${\bf a}_j:=(b_j,\epsilon_j),\,j\in \{1,2\}$ and
 ${\bf a}:=(b,\epsilon) =:{\bf a}_1\circ {\bf a}_2$,
then
$$
b=\frac{b_1 \overline{\epsilon}_2+b_2} {1+b_1\overline{b}_2\overline{\epsilon}_2}=
B_{(-b_2,1)}(b_1\overline{\epsilon}_2) ,\ \ \
\epsilon=\epsilon_1\frac{\epsilon_2+b_1\overline{b}_2} {1+\epsilon_2\overline{b}_1
b_2}=B_{(-b_1\overline{b}_2,\epsilon_1)} (\epsilon_2). \leqno( 1.3)
$$
 The neutral element of the group $(\mathbb B,\circ)$ is $e:=(0,1)\in \mathbb
 B$
and the inverse element of ${\bf a}=(b,\epsilon)\in \mathbb B$ is
${\bf a}^{-1}=(-b\epsilon,\overline \epsilon)$.

The integral of the  function  $f:\mathbb B\to \mathbb C$, with
respect to the left invariant Haar measure $m$ of the group $(\Bbb
B,\circ)$
 can be expressed as
  $$ \int_\mathbb B f({\bf a})\, dm({\bf a})=\frac
1{2\pi} \int_\mathbb I\int_\mathbb D \frac{f(b, e^{it})} {(1-|b|^2)^2}\, db_1 db_2
dt,\ \ \leqno( 1.4)
$$
where ${\bf a}=(b,e^{i t})=(b_1+ib_2,e^{it})\in \mathbb D\times\mathbb T$.

 It can be shown that
 this integral is invariant with respect to the left translation
${\bf a}\to {\bf a}_0\circ{\bf  a}$ and under the inverse transformation
 ${\bf a}\to {\bf a}^{-1}$, so this group is unimodular.

\vskip4mm

\centerline{\bf 2.1 Special discrete subsets  in $\Bbb B$ and their  sampling property }
\vskip5mm
The one parameter subgroups
$$
\Bbb B_1:=\{(r,1):r\in (-1,1)\},\ \ \Bbb B_2:=\{(0,\epsilon):\epsilon\in \Bbb T\} \leqno (2.1)
$$
generate $\Bbb B$, i. e.
$$
{\bf a}=(0,\epsilon_2)\circ (0,\epsilon_1)\circ (r,1)\circ (0,\overline \epsilon_1)\ \
({\bf a}=(r\epsilon_1,\epsilon_2),
\ r\in [0,1),\epsilon_1,\epsilon_2 \in \Bbb T). \leqno (2.2)
$$
$\Bbb B_1$ is the analogue of the group of dilation,
$\Bbb B_2$ is the analogue of the group of translation (see \cite{Sc09}).

The group operation $(r,1)=(r_1,1)\circ(r_2,1)$ in $\Bbb B_1$ can be expressed using the tangent hyperbolic and its inverse ($\text{ath}$)   in the following way 
$$
r=\frac{r_1+r_2}{1+r_1 r_2}=\text{th}(\text{ath}\ r_1+\text{ath}\ r_2)\ \ (r_1,r_2\in (-1,1)). \leqno (2.3)
$$

Let denote $r = \text{th} \,\alpha, r_i = \text{th} \,\alpha_i,\, i=1, 2 $. Then by
$$(r_1,1)\circ(r_2,1)=(\text{th} \,\alpha_1, 1)\circ (\text{th}\, \alpha_2, 1)= (\text{th} \,(\alpha_1+\alpha_2), 1),$$
it follows that  $(\Bbb B_1, \circ)$ is isomorphic to $(\Bbb R, +)$.
It is known that $(\Bbb Z, +)$ is a  subgroup of $(\Bbb R, +)$, then
$\overline{\Bbb B_1} =\{(\text{th}\, k, 1),\,\ k\in \Bbb Z \}$ is an
one parameter  subgroup of  $(\Bbb B_1, \circ)$ (see \cite{sobosc09}).

Let $a> 1$, denote by
$$
\Bbb B_3= \left\{(r_k, 1):\,\, r_k=\frac{a^k-a^{-k}}{a^k+a^{-k}}, \,
k\in \Bbb Z  \right\}.\leqno(2.4)
$$
It can be proved that $(\Bbb B_3, \circ)$ is another  subgroup of
$(\Bbb B, \circ)$, and $(r_k, 1)\circ (r_n, 1)=(r_{k+n}, 1)$. The
hyperbolic distance of the points $r_k,\ r_n$ has the following
property:
$$
\rho(r_k, r_n):=\frac{|r_k-r_n|}{|1-r_k\overline {r_n}|}=\left |\frac{\frac{a^k-a^{-k}}{a^k+a^{-k}}-\frac{a^n-a^{-n}}{a^n+a^{-n}}}{1-\frac{a^k-a^{-k}}{a^k+a^{-k}}\frac{a^n-a^
{-n}}{a^n+a^{-n}}}\right |= |r_{k-n}|. \leqno (2.5)
$$
Let $N(a,k), k\ge 1$ an increasing sequence of natural numbers, $N(a,0):=1$, and consider the following set of points 
 $z_{00}:=0$,
$$
{\cal{ A}}=\{ z_{k\ell}=r_k e^{i\frac{2\pi \ell}{N}},\,\,\,\, \ell = 0, 1, ..., N(a,k)-1,\,\,\,k=0,1, 2,... \} \leqno (2.6)
$$
and for a fixed $k\in \Bbb N$ let  the level $k$ be
$$
{\cal {A}}_k=\{ z_{k\ell}=r_k e^{i\frac{2\pi \ell}{N}},\, \ell \in \{ 0, 1, ..., N(a,k)-1 \} \, \}.\leqno (2.7)
$$
The points of ${\cal{A}}$  determine a similar,  decomposition to the Whitney cube decomposition of the unit disc   (see for ex. \cite{pa97} pp.80). There are two differences: 

1. The radius of the concentric circles  are connected to the Blaschke group operation, which is important from the point of view when we generate a multiresolution:  $(r_k, 1)\circ (r_n, 1)=(r_{k+n}, 1)$. This property   is  analogue with the property of  the dilatation  when we  generate  affine wavelet  multiresolution levels.  

2. The second:  we can choose $a$ and $N=N(a,k)$ such that ${\cal{A}}$ will be a set of  sampling in the Bergman space.  In the case of the Whitney cube decomposition  the concentric circles are defined by $r_k^*=1-\frac{1}{2^k}< r_k$ and we dived this in  $2^k$ equal parts, but the upper and lower density of this set it is not so easy to handle.

First we will study the following questions: for which choice  of $a$ and $N=N(a,k)$  

1. will be  ${\cal{A}}$ uniformly discrete,

2. will be  ${\cal{A}}$ an  $\epsilon$-net set for some $0<\epsilon < 1$,

3. will be ${\cal{A}}$ sampling sequence for Bergman spaces $A^p$?

To answer these questions let start with some basic definitions and results. For detailed exposition concerning these results see for example in \cite{dusc04}, \cite{scva}, \cite{duscvu}.

 Recall that if $z=x+iy\in\Bbb D$ then the normalized area measure is  $d A(z)=\frac{1}{\pi}d x d y$. For  $0<p<\infty$,
an  analytic   function $f:\Bbb D\to  \Bbb C$ belongs to the $A^p$ if 
$$\int_{\Bbb
D}|f(z)|^p\, d A(z)<\infty .\leqno (2.8)
$$
For $p=2$  the set $A^2$ is  a reproducing kernel Hilbert space, which is called the Bergman space. The reproducing kernel of $A^2$ is given by the formula 
$$
K(z.w)=\frac{1}{(1-\overline{w}z)^2}.\leqno (2.9)
$$
  The pseudohyperbolic metric is defined by
 $$
 \rho(z,y)=\left |\frac{y-z}{1-\overline{y}z} \right | \,\, (y,z \in \Bbb D).
 $$
A sequence of points $\Gamma=\{ z_k\}$ of points in the unit disc is uniformly discrete if
 $$
 \delta({\Gamma})=\inf_{j\ne k} \rho(z_j,z_k)=\delta>0.
 $$
 For $0<\epsilon<1$, a sequence of points $\Gamma=\{ z_k: k\in \Bbb N\}$ of points in the unit disc is said to be $\epsilon$-net if each point $z\in \Bbb D$
 has the property $\rho(z,z_k)<\epsilon$ for some $z_k$ in $\Gamma$. An equivalent statement is that 
 $$
 \Bbb D= \Cup_{k=1}^{\infty}\Delta (z_k, \epsilon),
 $$
 where $\Delta (z_k, \epsilon)$ denotes a pseudohyperbolic disc.
 
 A  sequence of points $\Gamma=\{ z_k: k\in \Bbb N\}$ of points in the unit disc is sampling sequence for $A^p$, where $0<p<\infty$, if there exist positive constants $A$ and $B$ such that 
 $$
 A ||f||^p \le \sum_{k=1}^{\infty} |f(z_k)|^p (1-|z_k|^2)^2 \le B ||f||^p, \quad 
f\in A^p. \leqno (2.10)
 $$
For $p=2$ this is   equivalent with
$$
 A ||f||^2 \le \sum_{k=1}^{\infty}\left| \langle f, \varphi_k \rangle \right|^2  \le B ||f||^2, \quad f \in A^2, \leqno (2.11)
 $$
 where $\varphi_k(z)=K(z,z_k)/\|K(z,z_k)\|$   are the normalized Bergman kernels with $w=z_k$ in (2.9). These kernel functions are not mutually orthogonal, so far no sequence of distinct points $z_k\in \Bbb D$ do the normalized kernel functions form an orthonormal basis. However, this last inequality shows that these functions will constitute a frame for $A^2$ if and only if $\Gamma=\{z_k: k\in \Bbb N\}$ is  a sampling set for $A^2$. A main difference between the Hardy space and the Bergman space is that there is no counterpart of sampling sequences in Hardy spaces. The Bergman spaces $A^p$ do have sampling sequences. but examples are not so easy to construct. Some explicit examples are due to Seip, Duren, Schuster, Horowitz, Luecking (see for ex in \cite{dusc04}). For our purpose we need a sampling sequence  connected to  the Blaschke  group operation.

%Closely related to the concept of sampling, and in some sense dual to it, is that of interpolation.
%A  sequence of points $\Gamma=\{ z_k: k\in \Bbb N\}$ of points in the unit disc is interpolation sequence for $A^p$ if the interpolation problem $f(z_k)=w_k$ for $k=1,2,...$ has a solution $f\in A^p$ whenever
%$$
%\sum_{k=1}^{\infty}(1-|z_k|^2)^2 |f(z_k)|^p< \infty.
%$$
%Every interpolation sequence is a zero-set in $A^p$, which means that there exists an $f\in A^p\setminus \{0\}$ which is zero on $\{ z_k: k\in \Bbb N\}$. 
 An $A^p$ sampling sequence is never an $A^p$ zero-set.  A  total characterization of sampling  sequences can be given with the uniformly discrete property and upper and lower density of the set see \cite{dusc04}, but these densities can be quite difficult to compute.  Duren, Schuster and Vukotic in \cite{duscvu} gave  for  sampling  sufficient conditions based on the pseudohyperbolic metric, that are relatively easy to verify, i.e: for $0<p<\infty$, if $\Gamma$ is a uniformly discrete $\epsilon$-net with 
$$
\epsilon<\frac{1}{1+\sqrt{\frac{2}{p}}}, 
$$
then is a sampling set for $A^p$. 

 Schuster and Varolin \cite{scva} improved these sufficient condition. They showed that every uniformly discrete $\epsilon$-net sequence with
$$
\epsilon<\sqrt{\frac{p}{p+2}} \leqno(2.12)
$$
is sampling set for $A^p$. 
This last sufficient condition will be used to answer our last questions.

{\bf Theorem 2.1} Let $a>1$,  and $(N=N(a,k), k\ge 1) $ a sequence of increasing natural numbers and  consider the  set of points
$$
{\cal {A}}=\{ z_{k\ell}=r_k e^{i\frac{2\pi \ell}{N}},\,\,\,\, \ell = 0, 1, ..., N(a,k)-1,\,\,\,k=0,1, 2,...\infty \}, \leqno 
$$
defined as before. Suppose that there exists %$N(a,k)(a^{2k}-a^{-2k})$ is a nondecreasing sequence, 
 $\alpha=\lim_{k\to \infty} N(a,k)a^{-2k}$.

1.
If $(N(a,k)a^{-2k}, k\ge 1) $ is increasing sequence and $\alpha$ is finite, then   ${\cal{A}}$ is uniformly discrete and the     separation constant satisfies
$$
\delta \ge \min \left \{ r_1, \frac{1}{\sqrt{1+\alpha^2}}\right \}.
$$

2.
If $(N(a,k)a^{-2k}, k\ge 1) $ is decreasing and $0<\alpha<\infty$, then there exists $\epsilon_0\in (0,1)$ for which the set $A$ is $\epsilon_0$-net  .

3. If $N(a,k)a^{-2k}=\alpha$, is constant for  $ k\ge 1$, $0<\alpha<\infty$ and 
$$
(a-a^{-1})^2+\pi^2\frac{a^2}{\alpha^2}<2p,
$$
then $A$ is a sampling set for $A^p$.

{\bf Proof}
1. We need to consider two types of situations. The pair of points lie on different circles, or they may lie on the same circler of radius $r_k$.  Suppose first that the points $z_{k\ell}, z_{mn}$ lie on two different circles of radius $r_k$ and $r_m$.
Then the generalized triangle inequality for the pseudohyperbolic metric (see \cite{dusc04} pp. 38) implies that
$$
\rho(z_{k\ell}, z_{mn})\ge \left| \frac{r_k-r_m}{1-r_kr_m}\right|\ge |r_{m-n}|\ge r_1>0.
$$
Next suppose that the pair of points lie on  the same circler of radius $r_k$, and $\ell \ne n$, then the least pseudohyperbolic distant is attained when $\ell=n+1$, then
$$
\rho(z_{k\ell}, z_{kn})\ge r_k \left|1- e^{\frac{2\pi i}{N}} \right|\left|1- r_k^2e^{\frac{2\pi i}{N}} \right|^{-1}
$$
$$
=2r_k \sin \frac{\pi}{N}\left[(1-r_k)^2+4r_k^2 \sin^2 \frac{\pi}{N}\right]^{-1/2}=
$$
$$
=\left\{ 1+[(1-r_k^2)/(2r_k\sin(\pi/N))]^2\right\}^{-1/2}.
$$
But $\sin(\pi/N)\ge (2/\pi)(\pi/N)=2/N$, so we deduce that 
$$\rho(z_{k\ell}, z_{kn})\ge\left\{ 1+[(1-r_k^2)N/(4r_k)]^2\right\}^{-1/2}.$$
%We keep in mind that $a^{-2k}N(a,k)$ is a nondecreasing sequence,
We observe that 
$$
(1-r_k^2)N/(4r_k)=\frac{N}{(a^{2k}-a^{-2k})}=N a^{-2k}[1/(1-a^{-4k})],
$$
and $\rho(z_{k\ell}, z_{kn})$ has a positive lower bound if  $\alpha=\lim_{k\to \infty}N(a,k)a^{-2k}<\infty$ and $\rho(z_{k\ell}, z_{kn})\ge \frac{1}{\sqrt{1+\alpha^2}}$.
Combining the two lower bounds we obtain the stated result for the separation constant.

2. For given $z=re^{i\theta}\in \Bbb D$ take $k$ and $j\in \{0,1,\cdots N(a,k)-1\}$ such that $r_k < r \le r_{k+1}$, $\theta \in [\frac{2\pi j}{N}, \frac{2\pi (j+1)}{N})$, $\theta_{kj}=\frac{2\pi j}{N}$, then 

$$
\frac{1}{1-\rho^2(z, z_{kj})}=\frac{(1-rr_k)^2+4rr_k\sin^2 \frac{\theta-\theta_{kj}}{2}}{(1-r^2)(1-r_k^2)}= 1 +\frac{(r-r_k)^2+4rr_k\sin^2 \frac{\theta-\theta_{kj}}{2}}{(1-r^2)(1-r_k^2)}\le
$$
$$
1 +\frac{(r-r_k)^2+4rr_k\frac{\pi^2}{N^2}}{(1-r^2)(1-r_k^2)}=1+\frac{(a-a^{-1})^2}{4}+\frac{(a^{2k+2}-a^{-2k-2})(a^{2k}-a^{-2k})}{4}\frac{\pi^2}{N^2}=:K(a,k).
$$
If  $(N(a,k)a^{-2k}, k\ge 1)$ is decreasing  and $\alpha=\lim_{k\to \infty} N(a,k)a^{-2k} \in (0,\infty)$, then the last term is upper bounded  by
$$
K:=1+\frac{(a-a^{-1})^2}{4}+\frac{a^2}{4\alpha^2}\pi^2.
$$
Then for 
$$
\epsilon_0=\sqrt{1-1/K},
$$ we have $\rho(z, z_{kj})< \epsilon_0.$

3.Using (2.12) we have that, if 
$(a-a^{-1})^2+\pi^2\frac{a^2}{\alpha^2}<2p$, then $\epsilon_0<\sqrt{\frac{p}{p+2}}$, consequently ${\cal{A}}$ is sampling sequence for $A^p$.

{\bf Remark} From this theorem we obtain that if ${\cal{A}}$ is a sampling set for $A^p$, then 
$$
(a-a^{-1})^2< 2p
$$
therefore $a$ must be in the interval  $(1, \frac{\sqrt{2p}+\sqrt{2p+4}}{2})$. Then we can always choose $N=N(a,k)$ big enough,  such that  the  the sampling condition will be satisfied. From the point of view of computations and to have on every circle the less possible numbers, for $p=2$ a convenient choice is  $a=2$, and $N(2,k)=2^{2k+\beta}$ for $k\ge 1$ with $\beta$ a fixed integer. Then $\alpha=2^{\beta}$, and the smallest value for $\beta$ for which the sampling condition is satisfied is  $\beta=3$,  then   on the $k$-th circle we will have  $N_1(2,k)=2^{2k+3}$ equidistant points corresponding to the roots of order $2^{2k+3}$ of the unity.   For $a=\sqrt{2}$ for sampling we need $N_1(\sqrt{2},k)=2^{k+2}$ points.

\centerline{\bf 2.2 The continuous voice transform on  Bergman space}

Our goal is to construct a multiresolution analysis based on the set ${\cal{A}}$ in the Bergman space.  First we will need a few results connected to the Bergman space \cite{dusc04} and the voice transform of the Bergman space.

The set $A^2=A^2(\Bbb D)$ with the scalar product
$$
\langle f,g\rangle:=\int_{\Bbb D}f(z)\overline{g(z)}
\, d A(z)\ \  \leqno (2.14)$$
 is a Hilbert space. An analytic  function in the unit disc of the form
$$
f(z):=\sum_{n=0}^\infty c_n z^n\ \ \ (z\in\Bbb D)
$$
 belongs to the set $A^2$ if and
only if the coefficients  satisfies
$$
\sum_{n=0}^\infty |c_n|^2 \frac 1{(n+1)}<\infty.
$$

The Bergman space   $A^2(\Bbb D)$ is a closed subspace of
$L^2(\Bbb D)$.

For each $z\in \Bbb D$ the point-evaluation map
$$
\tau_z:A^2(\Bbb D)\to \Bbb C,\quad \tau_z(f)=f(z)
$$
is a bounded linear functional on $A^2(\Bbb D)$. Each
function  $f\in A^2(\Bbb D)$ has the property
$$
|f(z)|\leqq \pi^{-1/2}\delta(z)^{-1}\|f\|_{A^2(\Bbb
D)} \quad (z\in \Bbb D),
$$
where $\delta(z)=dist(z, \Bbb T)$. From this it follows that the
 norm convergence in  $A^2(\Bbb D)$ implies the locally uniform
convergence on $\Bbb D$.
 Therefore,
by the Riesz Representation Theorem there is a unique element in
$A^2(\Bbb D)$, denoted by $K(., z)$, such that
$$
f(z)= \tau_z(f)=\langle f, K(., z)\rangle,\,\, (f\in
A^2(\Bbb D),  z\in \Bbb D).
$$
The function
$$
K: \Bbb D \times \Bbb D \to \Bbb C \quad \text{with}\quad K(., z)\in
A^2(\Bbb D)
$$
is called the Bergman kernel for $\Bbb D$. Taking $f(\xi)={K(\xi,
w)}$ for some $w\in \Bbb D$ we conclude that the kernel function has
the following symmetry property:
$$
{K(z,w)}=\frac1{\pi}\int_{\Bbb D}K(\xi,w)\overline{K(\xi, z)}d\xi_1
d\xi_2= \overline{K(w,z)}.
$$

 For any orthonormal basis $\{
\varphi_j, j=0,1,2,...\}$ for $A^2(\Bbb D)$ one has the
representation
$$
K(\xi,
z)=\sum_{j=1}^{\infty}\varphi_j(\xi)\overline{\varphi_j(z)},\quad
(\xi, z)\in \Bbb D\times \Bbb D,\leqno (2.15)
$$
with uniform convergence on compact subsets of $\Bbb D\times \Bbb
D$. The set of functions $$\{ \varphi_j(z)=\sqrt{(j+1)}z^j,\, z\in \Bbb D,\,
j=0,1,2,...\}$$ form an orthonormal basis in $A^2(\Bbb D)$,
consequently
$$
K(\xi, z)= \frac1{(1-\overline{z}\xi)^2}  .
$$
The explicit formula for the kernel function shows that
$$
f(z)= \frac1{\pi}\int_{\Bbb
D}f(\xi)\frac1{(1-\overline{\xi}z)^2}d\xi_1 d\xi_2 \quad (f\in
A^2(\Bbb D),\,z\in \Bbb D). \leqno (2.11)
$$
Applying this formula in particular for
$f(z)=(1-\xi\overline{z})^{-2}$ for fixed $z$ in the disc we obtain
that
$$
\|K(y,.)\|^2_2=\frac1\pi \int_{\Bbb
D}\frac1{|1-\overline{\xi}z|^4d\xi_1d\xi_2}=\frac1{(1-|z|^2)^2}=K(z,z)>0.
$$

The voice transform on Bergman space   is induced by a  unitary
representation of the Blaschke group on the  Bergman
space. Results connected  the voice transform on Bergman space were published in \cite{pasc09}.

Let consider the following set of functions
$$
F_{\bf a}(z):=\frac{\sqrt{\epsilon(1-|b|^2)}}{1-\bar b z}\ \ \
(a=(b,\epsilon)\in\Bbb D, z\in\overline {\Bbb D}).\leqno 2.10
$$
This function induces a unitary
representation on the space $A^2$. Namely let define
$$
U_{\bf a}f:=[F_{a^{-1}}]^2 f\circ B^{-1}_a\ \ \ ({\bf a}\in\Bbb B, f\in A^2),
$$
then this is a representation, i.e.:
$$
i)  \ \ U_{x\cdot y}=U_x\circ U_y\ \ (x,y\in \Bbb B),$$
 $$ii) \ \ \Bbb B\ni x\to U_xf\in A^2\ \ \textrm{is continuous for all} \ f\in A^2. \leqno
$$
And more  $U_{\bf a} ({\bf a}\in\Bbb B)$ is an unitary, irreducible square integrable 
representation of the group  $\Bbb B$ on the Hilbert space $A^2$.

The {\it voice transform } of $f\in A^2$ generated by the
representation $U$ and by the parameter  $g\in A^2$ is the
(complex-valued)  function on $\Bbb B$ defined by
$$
{(V_g f)(x):=\langle f,U_x g\rangle\ \ \ (x\in \Bbb B, f, g  \in
A^2).}\leqno 
$$

In  \cite{pasc09} a direct proof of the analogue of Plancherer formula  was given, i.e.:
$$
[V_{\rho_1} f,V_{\rho_2} g]= 4\pi \langle \rho_1, \rho_2\rangle\
\langle f, g\rangle\ \ (f,g,\rho_1, \rho_2\in A^2(\Bbb
D)),
$$
where
$$
[F, G]:=\int_{\Bbb B}F(a)\overline{G(a)}\, dm(a),
$$
$m$ is the Haar measure of the group $\Bbb B$.

This transform is in same relation with the Blaschke group and the Bergman space as
 the affine wavelet transform with the affine group and the $L^2(\Bbb R)$ (see \cite{Sc09}).  Indeed let consider $(G, \circ)$ equal to the  the affine group, where
$$
G=\{\ell_{(a,b)}:\Bbb R\to\Bbb R:  (a,b)\in \Bbb R^*\times \Bbb R \},$$ $$
\ell_{(a,b)}(x)=ax+b, \,\,\Bbb R^*:=\Bbb R \setminus \{0\},\,\,
\ell_1 \circ \ell_2(x)=\ell_1(\ell_2(x))=a_1a_2x+a_1b_2+b_1.$$
  The representation of the affine group $G$ on $L^2(\Bbb R)$ is given by
$$
U_{(a,b)}f(x)=|a|^{-1/2}f(a^{-1}x-b),
$$
where $a$ is the dilatation parameter, and $b$ the translation parameter.

 The continuous affine wavelet transform is a voice transform generated
by this representation of the  affine group:
$$
W_{\psi}f(a,b)=|a|^{-1/2}\int_{\Bbb R}f(t)\overline{\psi(a^{-1}t-b)}dt
=\langle f,U_{(a,b)}\psi\rangle,\quad f, \psi \in L^2(\Bbb R).
$$
There is a rich bibliography of the affine wavelet theory (see for example \cite{ch92, da88, hewa89, ki, me90, ma89, scwa95}.% [ 7, 9, 17, 19, 20,  21, 34] etc.).
 One important question is the construction of the discrete version, i.e.,  to find $\psi$ so that the discrete translates and dilates
$$
\psi_{n,k}=2^{-n/2}\psi(2^{-n}x-k)
$$
form a (orthonormal) basis in $L^2(\Bbb R)$ and generate a multiresolution (see [7 , 9, 17, 19, 20, 21] etc.).
 The construction of the discrete version is connected to  the discrete subgroup of the affine group, generated by the following set:
 $$
G_{n,k}=\{\ell_{(2^{-n},-k)}:\Bbb R\to\Bbb R:  n\in \Bbb Z,\, k\in \Bbb Z\}.
$$
   The discretization of the voice transform can be achieved using the unified approach of the atomic decomposition elaborated by Feichtinger and Gr\"ochenig \cite{fegr89-1}. This general description can be applied when the integrability condition of the voice transform is satisfied. In a recent paper  \cite{pa12} it is  shown  that, the integrability condition in the Bergman space it is not satisfied. This motivates to find discrete multiresolution decomposition in Bergman spaces?

 \vskip5mm
\centerline{\bf 2.3 Multiresolution analysis in  the Bergman space}
\vskip5mm

 We start with the general definition of the  affine wavelet multiresolution analysis in $L^2(\Bbb R)$.

{\bf Definition 2.2.1.} {\it Let $V_j, \,\, j\in \Bbb Z$ be a sequence of subspaces of $L^2(\Bbb R)$. The collections of spaces $\{V_j, \,\, j\in \Bbb Z\}$ is called a multiresolution  analysis with scaling function $\phi$ if the following conditions hold:

1. (nested) $V_j\subset V_{j+1}$

2. (density) $\overline{\cup V_j}= L^2(R)$

3. (separation) $\cap V_j=\{ 0\}$

4. (basis) The function $\phi$ belongs to $V_0$ and the set $\{ 2^{n/2}\phi(2^nx-k), \,\, k\in \Bbb Z\}$ is a (orthonormal) bases in $V_n$.}

In multiresolution analysis, one decomposes a function space in several resolution levels and the idea is to represent the functions from the function space by a low resolution approximation and adding to it the successive details that lift it to resolution levels of increasing detail.

Wavelet analysis couples the multiresolution idea with a special choice of bases for the different resolution spaces and for the wavelet spaces that represent the difference between successive resolution spaces.  If $V_n$ are the resolution spaces $V_0\subset V_1 \subset....\subset V_n...$, then the wavelet spaces $W_n$ are defined by the  equality $W_n\bigoplus V_n=V_{n+1}$.

In the construction of affine wavelet multiresolutions the dilatation is used to obtain a higher level resolution ($f(x)\in V_n \Leftrightarrow f(2x)\in V_{n+1}$) and applying the   translation we remain on the same level of resolution. This field has also a rich bibliography (see for example \cite{boscsz99, bugo99, da88, ch92, dupa96, evca98, hewa89, ki, me90, ma89, wapa96, wapa98}).
Using the subgroup $\Bbb B_3$ of the Blaschke group, a discrete subgroup of $\Bbb B_2$ and the representation we give a similar construction of the affine wavelet multiresolution in the Bergman space. 
 To show the analogy with the affine wavelet multiresolution  we first represent the levels $V_n$  by nonorthogonal bases and then we construct an orthonormal bases in $V_n$ and  give also  an orthogonal basis   in $W_n$ which is orthogonal to $V_n$. We will show that  in the case of this discretization the analogue of the  Malmquist-Takenaka systems for Bergman space, will span the resolution spaces and the density property  will be fulfilled, i.e., $\overline{\bigcup_{k=1}^{\infty}V_k}=A^2$ in norm. Similar multiresolutin results for the Hardy space ${\cal{H}}^2(\Bbb T)$ were obtained by the author in \cite{pa11}.

We show that  the  projection $P_nf$  on the $n$-th resolution level is an interpolation operator in the unit disc until the $n$-th level, which converges in $A^2$ norm to $f$.

Let consider
$a>1$, denote by
 $r_k=\frac{a^k-a^{-k}}{a^k+a^{-k}}, \,
k\in \Bbb N $,  $N=N(a,k)$ a sequence of  natural numbers such that 
 $\alpha= N(a,k)a^{-2k}$ satisfies
$$0<\alpha<\infty, \quad
(a-a^{-1})^2+\pi^2\frac{a^2}{\alpha^2}<4.
$$
Let us consider the  set of points
$$
{\cal{A}}=\{ z_{k\ell}=r_k e^{i\frac{2\pi \ell}{N}},\,\,\,\, \ell = 0, 1, ..., N-1,\,\,\,k=0,1, 2,... \}, \leqno 
$$
and for a fixed $k\in \Bbb N$ let  the level $k$ be
$$
{\cal{A}}_k=\{ z_{k\ell}=r_k e^{i\frac{2\pi \ell}{N}},\, \ell \in \{ 0, 1, ..., N-1 \} \, \}.
\leqno (2.17)
$$
Due to Theorem 2.1  ${\cal{A}}$ is a sampling set for $A^2$. This implies that the set of normalized kernels $(\varphi_{kl}(z)=\frac{(1-r_k^2)}{(1-\overline{z_{k\ell}}z)^2},\, k=0,1,\cdots, \, \ell=0,1,\cdots N-1)$ will constitute a frame system for $A^2$. From the frame theory \cite{gr01} or from atomic decomposition results (see  Theorem 3 of \cite{zh95}) follows that every function $f$ from $A^2$ can be represented
$$
f(z)=\sum_{(k,\ell)}c_{k\ell}\varphi_{kl}(z)
$$
for some $\{c_{k\ell}\}\in \ell^2$, with the series converging in $A^2$ norm. The determination of the coefficients it is related to the construction of the inverse frame operator (see \cite{gr01}), which is not an easy task in general. This is the reason why we try to construct other approximation process for $f\in A^2$ such that the determination of the coefficients follow an exactly defined algorithmic scheme.
 
Let us consider the  function $\varphi_{00}=1$ and let $V_0=\{ c,\,\,
c\in \Bbb C \}$. 
 
Let us consider the nonorthogonal hyperbolic wavelets at the first
level
$$ \varphi_{1,\ell}(z)=(U_{(z_{1\ell},1)^{-1}}p_0)(z)=\frac{{(1-r_1^2)}}{(1-\overline{z_{1\ell}}z)^2},\quad \ell=0,1,\cdots N(a,1).\leqno(2.18)$$
 They can be obtained from $\varphi_{1,0}$ using the analogue of translation operator which in the unit disc is a multiplication by a unimodular complex number, and from $\varphi_{0,0}$  using first the  representation operator $U_{(r_1,1)^{-1}}$ followed by the translation operator:
  $$\varphi_{1,\ell}(z)=\varphi_{1,0}(ze^{-\frac{2\pi i \ell}{N(a,1)})})=(U_{(r_1,1)^{-1}}\varphi_{0,0})(ze^{-\frac{2\pi i\ell}{N(a,1)})}).\leqno(2.19)$$
   Let us define the first resolution level as follows
$$
V_1=\{ f: \Bbb D\to \Bbb C,\,\, f(z)=c_{0,0}\varphi_{0,0}+\sum_{\ell=0}^{N(a,1)-1}c_{1,\ell}\varphi_{1,\ell},\,\, c_{00}, c_{1,\ell}\in \Bbb C, \, \ell=0,1,\cdots N(a,1)-1. \}.\leqno (2.20)
$$
Let us consider the nonorthogonal wavelets on the $n$-th level
$$ \varphi_{n,\ell}(z)=(U_{(z_{n\ell},1)^{-1}}p_0)(z)=\frac{{(1-r_n^2)}}{(1-\overline{z_{n\ell}}z)^2},\quad \ell=0,1,..., N(a,n)-1,\leqno (2.21)$$
 which can be obtained from $\varphi_{n,0}$ using the  translation  operator, and  from $\varphi_{0,0}$ using the representation $U_{((r_{n-1},1)\circ (r_1,1))^{-1}}$, and the  translations  $$\varphi_{n,\ell}(z)=(U_{((r_{n-1},1)\circ (r_1,1))^{-1}}p_0)(ze^{-i\frac{2\pi \ell}{N(a,n)})}).\leqno (2.22)$$
Let us define the $n$-th resolution level by
$$
V_n=\{ f: \Bbb D\to \Bbb C,\,\, f(z)=\sum_{k=0}^n\sum_{\ell=0}^{N(a,k)-1}c_{k,\ell}\varphi_{k,\ell},\,\,  c_{k,\ell}\in \Bbb C\  \}.\leqno (2.23)$$
The closed subset $V_n$ is spanned by
$$
\{ \varphi_{k,\ell},\,\,\ell=0,1,...,N(a,k)-1,\,\, k=0,...,n  \}.\leqno (2.24)
$$
Continuing this procedure we obtain a sequence of closed, nested subspaces of $A^2$ for $z\in\Bbb D$
$$
V_0\subset V_1 \subset V_2 \subset .....V_n \subset ....\, A^2.\leqno (2.25)
$$

Due to Theorem 2.1 the  normalized kernels $\{\varphi_{kl}(z)=\frac{(1-r_k^2)}{(1-\overline{z_{k\ell}}z)^2},\, k=0,1,\cdots, \, \ell=0,1,\cdots N(a,k)-1\}$ form a frame system for $A^2$   which implies that  this is a   complete and
 a closed set in norm, i.e.,
$$
\overline{\bigcup_{n\in \Bbb N} V_n}=A^2, \leqno (2.26)
$$
consequently the density property it is satisfied.

For $a=2$ and $N(2,n)=2^{2n+3}$, 
if a function $f\in V_n$, then $U_{(r_1,1)^{-1}}f \in V_{n+1}$. For this it is enough to show that
$$
U_{(r_1,1)^{-1}}(\varphi_{k,\ell})(z)=U_{(r_1,1)^{-1}}[(U_{(r_{k},1)^{-1}}p_0)](ze^{-i\frac{2\pi \ell}{2^{2k+3}})})=$$ $$
[(U_{(r_{k+1},1)^{-1}}p_0)](ze^{-i\frac{2\pi 4\ell}{2^{2(k+1)+3}}})\in V_{n+1},\quad k=1,..,n,\,\, \ell=1,...2^{2k+3}-1. \leqno (2.27)
$$
From now on for simplicity we will deal with this case, but in the general case  we can always choose $N(a,k)$ such that the previous condition to be valid. 

Since the set ${\cal{A}}$ is a sampling set it follows that is a set of uniqueness for $A^2$, which means that every function $f\in A^2$ is uniquely determined by the values $\{f(z_{k\ell})\}$. In the paper of Kehe Zhu \cite{zh97} described  in general how can be recaptured a function from a Hilbert space when  the values of the function on a set of uniqueness are known  and developed in details this process in the Hardy space. At the beginning we will follow the steps of the recapturation process but we will combine this with the multiresolutin analysis. 
The set
$$
\left\{ \frac{1}{(1-\overline{z_{k\ell}}z)^2},
\,\,\ell=0,1,...,2^{2k+3}-1,\,\, k=0,1,...,n.  \right\}\leqno (2.28)
$$
is a nonorthogonal basis in $V_n$.

Using Gram-Schmidt
orthogonalization process they can be orthogonalized. Denote by $\psi_{k,\ell}$ the resulting functions. They can be seen as the analogue of  the Malmquist -Takenaka system in the Hardy space. This functions can be obtained using the following two methods.
The first  arises from the orthogonalization procedure. To describe this let reindex the points of the set ${\cal{A}}$ as follows $a_1=z_{00}, a_2=z_{10}, a_3=z_{11}, a_4=z_{12}, ..., a_{33}=z_{1,31}, a_{34}=z_{2,0}, \cdots a_{m}=z_{k\ell}\cdots , k=0,1,...., \ell=0,1,...,2^{2k+3}-1$, and denote by $K(z,z_{k\ell})=\frac{1}{(1-\overline{z_{k\ell}}z)^2}:=K(z,a_m)$
\begin{equation}
\begin{cases}
\phi_{00}(z)=\phi(a_1,z)=\frac{K(z,a_1)}{\|K(.,a_1)\|},\\
\phi_{k\ell}(z)=\phi(a_1,a_2,...,a_m,z)=K(z, a_m)-\sum_{i=1}^{m-1}\phi({a_1,a_2,...,a_i,z})\frac{\langle K(.,a_m),\phi(a_1,a_2,...,a_i,.)\rangle }{\|\psi(a_1,a_2,...,a_i,.) \|^2},\\ a_m=z_{k\ell},\,  m\ge 2.
\end{cases}
\end{equation}
Thus the normalized functions $\{ \psi_{k\ell}(z)=\frac{\phi_{k\ell}(z)}{\|\phi_{k\ell} \|}, \, k=1,2,\cdots , \ell=0,1,\cdots 2^{2k+3}\}$ became  an orthonormal system.
Applying similar construction in Hardy space we  get in this way the Malmquist-Takenaka system. They can be written in a nice closed form using the Blaschke products. Unfortunately in our situation this is not the case and the properties of the system can be not seen from the previous construction.

Another approach is given by Zhu in \cite{zh97}. He  shows that the result of the Gram-Schmidt process  are connected to some reproducing kernels and the contractive divisors on Bergman spaces. Let denote $A_m=\{ a_1, a_2,\cdots a_m \}$. Let $H_{A_m}$ the subspace of $A^2$ consisting of all functions in $A^2$ which vanish on $A_m$. $H_{A_m}$ is a closed subspace of $A^2$ and denote by $K_{A_m}$ the reproducing kernel of $H_{A_m}$. These reproducing kernels satisfies the following recursion formula 
$$
K_{A_{m+1}}(z,w)=K_{A_{m}}(z,w)-\frac{K_{A_{m}}(z,a_{m+1})K_{A_{m}}(a_{m+1},w)}{K_{A_{m}}(a_{m+1},a_{m+1})}, m\ge 0,$$ $$ K_{A_{0}}:=K(z,w)=\frac{1}{(1-\overline{w}{z})^2}. \leqno(2.29)
$$
The result of the Gram-Schmidt process can be expressed as 
$$
\frac{K(z,a_1)}{\sqrt{K(a_1,a_1)}},\, \frac{K_{A_1}(z,a_2)}{\sqrt{K_{A_1}(a_2,a_2)}},\,\cdots \frac{K_{A_{m-1}}(z,a_{m})}{\sqrt{K_{A_{m-1}}(a_{m},a_{m })}},\cdots.
$$
Then 
$$
\psi_{k\ell}(z)=\frac{K_{A_{m-1}}(z,a_{m})}{\sqrt{K_{A_{m-1}}(a_{m},a_{m})}} \leqno(2.30)
$$
and is the solution of the following  problem
$$
\sup\{ Re f(a_{m}): f\in H_{A_{m-1}}, \|f\|\le 1\}.
$$
This extremal functions in the context of the Bergman spaces have been studied extensively in recent years by Hedenmalm \cite{he91}.
The main result in \cite{he91} is that the function
$$
\frac{K_{A_{m-1}}(z,a_{m})}{\sqrt{K_{A_{m-1}}(a_{m},a_{m})}}
$$
is a contractive divisor on the Bergman space, vanishes on $A_{m-1}$, and if ${\cal{A}}$ is not a zero set for $A^2$, as is in our case, the functions converge  to $0$ as $m\to \infty$.  In Hardy space the partial products of a Blaschke product corresponding to a nonzero set  own all these nice properties.

From the Gram-Schmidt orthogonalization process it follows that
$$
V_n= span \{ \psi_{k, \ell},\,\,\ell=0,1,..., 2^{2k+3}-1, \,\, k=\overline{0, n} \}.\leqno (2. 31)
$$

The wavelet space $W_n$ is the orthogonal complement of $V_n$ in $V_{n+1}$.
We will prove that
$$
W_n=span \{ \psi_{n+1, \ell},\,\,  \,\, \ell=0,1,..., 2^{2n+5}-1 \}.\leqno (2.32)
$$

Indeed, every function $f\in A^2$ can be recovered using the,
$$
f(z)=\int_{\Bbb D}\frac {f(w)}{(1-\overline{w}z)^2}dA(w). \leqno (2.33)
$$
If $f\in V_n$, one has
$f(z)=\sum_{k=0}^n\sum_{\ell=0}^{2^{2k+3}-1}c_{k,\ell}\varphi_{k,\ell}
 \subset A^2$, then using (2.30) we obtain that
$$
\langle  \psi_{n+1 j}, f\rangle=\sum_{k=0}^n\sum_{\ell=0}^{2^{2k+3}-1}c_{k,\ell}\langle  \psi_{n+1, j},\varphi_{k,\ell}\rangle=$$ $$ \sum_{k=0}^n\sum_{\ell=0}^{2^{2k+3}-1}c_{k,\ell}{(1-r_k^2)} \psi_{n+1, \ell}(z_{k\ell})=0,\,\, j=0,1,...2^{2n+5}-1.
$$
We have proved that for $f\in V_n$
$$
\langle f, \psi_{n+1, j}\rangle=0, \leqno (2. 34)
$$
which is equivalent with
$$
\psi_{n+1, j}\perp V_n, \quad (j=0,1,..., 2^{2n+5}-1).\leqno (2. 35)
$$
From
$$
V_{n+1}=V_n\bigoplus span \{ \varphi_{n+1, j},\,\, j=0,1,..., 2^{2n+5}-1\} \leqno (2.36)
$$
it follows that $W_n$ is an $2^{2(n+1)+3}$ dimensional space and
$$
W_n=span \{ \psi_{n+1, \ell},\,\,  \,\, \ell=0,1,..., 2^{2n+5}-1 \}. \leqno (2.37)
$$
{\bf Summary.}  We have generated a
multiresolution in $A^2$ and we have constructed a rational orthogonal wavelet system which generates the levels of the  multiresolution.
%\newpage

\vskip5mm \centerline{{\bf 2.4 The projection operator corresponding
to the $n$-th resolution level  }}

\vskip5mm Let us consider the
orthogonal projection operator of an arbitrary function $f\in
A^2$ on the subspace $V_n$ given by
$$
P_nf(z)=\sum_{k=0}^n\sum_{\ell=0}^{2^{2k+3}-1}\langle f, \psi_{k,\ell}\rangle \psi_{k,\ell}(z). \leqno (2.38)
$$
This operator is called the projection of $f$ at scale or resolution level $n$.
\vskip3mm
{\bf Theorem 2.2 } {\it For $f\in A^2$ the projection operator $P_nf$ is an
interpolation operator on the points $z_{k\ell}=r_k e^{i\frac{2\pi \ell}{2^{2k+3}}},\,\,(\ell=0,...., 2^{2k+3}-1,\,\,\ k=0,...,n)
$, is norm convergent in $A^2$ to $f$, i.e.
$$
\| f-P_nf \|\to 0, \quad n\to \infty, 
$$
 uniformly  convergent  inside the unit disc on every compact subset, and is the solution of minimal norm interpolation problem.}
\vskip3mm
{\bf Proof}  Let consider $N=1+2^5+\cdots + 2^{2n+3}$ and the corresponding  kernel function of the projection
operator
$$
{{\bf K}_N(z, \xi)}=\sum_{k=0}^n\sum_{\ell=0}^{2^{2k+3}-1} \overline{\psi_{k,\ell}(\xi)} \psi_{k,\ell}(z)=$$
$$\sum_{m=1}^N\frac{K_{A_{m-1}}(z,a_{m})}{\sqrt{K_{A_{m-1}}(a_{m},a_{m})}}\overline{\left ( \frac{K_{A_{m-1}}(\xi,a_{m})}{\sqrt{K_{A_{m-1}}(a_{m},a_{m})}}\right )}= \sum_{m=1}^N\frac{K_{A_{m-1}}(z,a_{m})K_{A_{m-1}}(a_{m}, \xi)}{{K_{A_{m-1}}(a_{m},a_{m})}}.\leqno (2.39)
$$
From the recursion relation (2.29) it follows that
$$
{{\bf K}_N(z, \xi)}=\sum_{m=1}^{N}(K_{A_{m-1}}(z,\xi)-K_{A_{m}}(z,\xi))=K(z,\xi)-K_{A_{N}}(z, \xi) \leqno(2.40)
$$

From this relation it follows that the values of the kernel-function in the points $z_{k\ell},\,\, (\ell=0,...., 2^{2k+3}-1,\,\,\ k=0,...,n)$ are equal to
$$
K({z_{kl}}, \xi)=\frac{1}{(1-z_{k\ell}\overline{\xi})^2}. \leqno (2.41)
$$
Every function $f\in A^2$ can be recovered using the Bergman projection
$$
f(z)=\int_{\Bbb D}\frac {f(w)}{(1-\overline{w}z)^2}d A(w)
$$
Therefore
$$
P_nf(z_{k\ell})=\int_{\Bbb D}\frac{f(w)}{(1-\overline{w}z_{mj})^2}d A(w)=f(z_{mj})
\quad (j=0,...., 2^{2m+3}-1,\,\,\ m=0,...,n). \leqno (2.42)
$$
 We obtain that $P_nf$ is  interpolation operator for every $f\in A^2$ on the set $\cup_{m=0}^n{\cal{A}}_m$.

Because of 2.26 and 2.31 $\{ \psi_{k, \ell},\,\, k=\overline{0, \infty}, \,\,
\ell=0,1,..., 2^{2k+3}-1 \}$ is a closed set in the Hilbert space
$A^2$, we have that
that $\| f-P_nf\|\to 0$ as $n\to \infty$. Since convergence in
$A^2$ norm implies uniform convergence  on every compact subset inside the unit disc, we
conclude that $ P_nf(z)\to  f(z)$ uniformly on every compact subset
of the unit disc. From Theorem 5.3.1 of \cite{pa97} there exists
 a unique $\hat{f}_n \in V_n$ with minimal norm such that
$$
\hat{f}_n(z_{mj})=f(z_{mj}),\quad (j=0,...., 2^{2m+3}-1,\,\,\
m=0,...,n), \leqno (2.43)
$$
 $\hat{f}_n$ is uniquely determined by the interpolation
conditions
and is equal to the orthogonal projection of $f$ on $V_n$, 
thus
 $\hat{f}_n(z)=P_nf(z)$.

 \vskip3mm 
 \centerline{\bf 2.5 Reconstruction algorithm }
\vskip3mm

In what follows we propose a computational scheme for the best
approximant in the wavelet base  $\{
\psi_{k,\ell},\,\,\ell=0,1,...,2^{2k+3}-1,\,\, k=0,...,n  \}$.

The projection of $f\in A^2$ onto $V_{n+1}$  can be  written in the following way:
$$
P_{n+1}f=P_nf+ Q_nf, \leqno (2.44)
$$
where
$$
Q_nf(z):=\sum_{\ell=0}^{2^{2n+5}-1}\langle f, \psi_{n+1,\ell}\rangle \psi_{n+1,\ell}(z). \leqno (2.45)
$$
This operator has the following properties
$$
Q_nf(z_{k\ell})=0, \quad k=1,...,n,\,\, \ell=0,1,..., 2^{2n+3}-1. \leqno (2.46)
$$
 Consequently $P_n$
contains information on low resolution, i.e., until the level ${\cal{A}}_n$, and $Q_n$ is the high resolution part.
After $n$ steps
$$
P_{n+1}f=P_1f+ \sum_{k=1}^n Q_nf. \leqno (2.47)
$$
Thus
$$
V_{n+1}=V_0\bigoplus W_0\bigoplus W_1\bigoplus...\bigoplus W_n. \leqno (2.48)
$$

The set of coefficients of the best approximant $P_nf$
$$\{ b_{k\ell}=\langle f, \psi_{k,\ell}\rangle, \, \,\ell=0.1,...2^{2k+3}-1\,\,\ k=0,1,..., n \} \leqno (2.49)$$
is the (discrete) hyperbolic wavelet transform of the function $f\in A^2$.
 Thus it is important to have an efficient algorithm for the computation of the coefficients.

The coefficients of the projection operator $P_nf$ can be computed
if we know the values of the functions on $\bigcup_{k=0}^n {\cal{A}}_k$. For
this reason  we express first the function $\psi_{k,\ell}$ using the
bases $(\varphi_{k',\ell'}\,\,\ell'=0,1,...2^{2k'+3}-1,\,\,
k'=0,...,k)$, i.e. we write the partial fraction decomposition of
$\psi_{k\ell}$ :
$$
\psi_{k,\ell}=\sum_{k'=0}^{k-1}\sum_{\ell'=0}^{2^{2k'+3}-1}c_{k',\ell'}\frac1{(1-\overline{z_{k'\ell'}}\xi)^2}+
\sum_{j=0}^{\ell}c_{k,j}\frac1{(1-\overline{z_{kj}}\xi)^2}. \leqno (2.50)
$$
Using the orthogonality of the functions
$(\psi_{k',\ell'}\,\,\ell'=0,1,...2^{2k'+3}-1,\,\, k'=0,...,k )$ and the reconstruction formula
$$
\delta_{kn}\delta_{\ell m}=\langle \psi_{nm}, \psi_{k\ell}\rangle=
\sum_{k'=0}^{k-1}\sum_{\ell'=0}^{2^{2k'+3}-1}\overline{c_{k',\ell'}}\psi_{n,m}(z_{k'\ell'})+
\sum_{j=0}^{\ell}\overline{c_{k,j}}\psi_{n,m}(z_{kj}), \leqno (2.51)
$$
$$
 (m=0,1,...2^{2n+3}-1,\,\,n=0,...,k).
$$
If we order these equalities so that  we write first the relations (2.51) for $n=k$ and  $m=\ell, \ell-1,...,0$ respectively, then for $n=k-1$ and $m=2^{2(k-1)+3}-1,\ 2^{2(k-1)+3}-2,...,0 $, etc., this is equivalent to
$$\left(\begin{matrix} 1\\0\\0\\.\\.\\.\\0\end{matrix}\right)
          =\left(\begin{matrix}%
\psi_{k,\ell}(z_{k,\ell}) & 0 & 0 & \dots & 0 \\
\psi_{k,\ell-1}(z_{k,\ell}) & \psi_{k,\ell-1}(z_{k,\ell-1}) & 0 & \dots & 0 \\
\psi_{k,\ell-2}(z_{k,\ell})& \psi_{k,\ell-2}(z_{k,\ell-1}& 0 & \dots & 0\\
\vdots & & &\vdots \\
\psi_{00}(z_{k,\ell}) &\psi_{00}(z_{k,\ell-1}) & \psi_{00}(z_{k,\ell-2}) & \dots & \psi_{00}(z_{00}) \\
\end{matrix}\right)\left(\begin{matrix}\overline{c_{k,\ell}}\\ \overline{c_{k,\ell-1}}\\\overline{c_{k,\ell-2}}\\ \vdots\\ \overline{c_{00}}\end{matrix}\right). \leqno (2.52)$$

This system has a unique solution $(\overline{c_{k,\ell}}, \overline{c_{k,\ell-1}},\overline{c_{k,\ell-2}}, ..., \overline{c_{00}})^T$. If we determine this vector, then we can compute the exact value of $\langle f, \psi_{k,\ell}\rangle$ knowing the values of $f$ on the set $\bigcup_{k=0}^n {\cal{A}}_k$.

Indeed, using again the partial fraction decomposition of $\psi_{k,\ell}$ and the reconstruction formula formula we get that
$$
\langle f,\psi_{k,\ell}\rangle=\sum_{k'=0}^{k-1}\sum_{\ell'=0}^{2^{2k'}-1}\overline{c_{k',\ell'}}\langle f(\xi)\frac1{(1-{z_{k'\ell'}}\overline{\xi})^2}\rangle+
\sum_{j=0}^{\ell}\overline{c_{k,j}}\langle f(\xi),\frac1{(1-{z_{kj}}\overline{\xi})^2}\rangle=
$$
$$
=\sum_{k'=0}^{k-1}\sum_{\ell'=0}^{2^{2k'}-1}\overline{c_{k',\ell'}} f(z_{k',\ell'})+
\sum_{j=0}^{\ell}\overline{c_{k,j}}f(z_{k,j}). \leqno (2.53)
$$
{\bf Summary.} Measuring the values of the function $f$ in the points of the set
 ${\cal{A}}=\bigcup_{k=0}^{n}{\cal{ A}}_k \subset \Bbb D$ we can write the operator $(P_nf, \, n\in \Bbb N)$ which
  is convergent in $A^2$ norm to $f$,  is the minimal norm  interpolation operator on the
  set the $\bigcup_{k=0}^n{\cal{ A}}_k$ and $ P_nf(z)\to  f(z)$
  uniformly on every compact subset of the unit disc.

\end{document}